\documentclass[12pt]{article}

\usepackage{a4wide}
\usepackage{amssymb,amsmath}
\usepackage{version}
\excludeversion{temp}

\allowdisplaybreaks[1]

\newcommand{\titre}{Quantised coordinate rings of semisimple groups \\are 
unique factorisation domains }

\newenvironment{proof}{\begin{trivlist}\item[]{\it
Proof.}}{\hfill$\square$\end{trivlist}}

\newtheorem{theorem}{Theorem}[section]
\newtheorem{corollary}[theorem]{Corollary}

\newtheorem{lemma}[theorem]{Lemma}
\newtheorem{proposition}[theorem]{Proposition}

\newcommand{\spec}{{\rm Spec}}

\newcommand{\g}{\mathfrak{g}}
\newcommand{\bk}{\mathfrak{b}}
\newcommand{\comp}{\mathbb{C}}

\newcommand{\vl}{V(\lambda)}
\newcommand{\dime}{\mathrm{dim}}

\newcommand{\ideal}[1]{\langle {#1}\rangle}

\def\ch{{\mathcal H}}

\begin{document}

\title{\titre}
\author{S Launois and  T H Lenagan
\thanks{This research was supported by a
 Marie Curie Intra-European Fellowship within the $6^{\mbox{th}}$
 European Community Framework Programme and by Leverhulme Research Interchange Grant F/00158/X }\;
}
\date{}

\maketitle


\begin{abstract}
We show that the quantum coordinate ring of a semisimple group is a unique
factorisation domain in the sense of Chatters and Jordan in the case where the
deformation parameter $q$ is a transcendental element.
\end{abstract}



\section*{Introduction}\label{sec:intro}

Throughout this paper, $\comp$ denotes the field of
complex numbers, $\comp^*:=\comp \setminus \{0\} $ and $q \in \comp^*$
is transcendental.

The notion of a noncommutative noetherian unique factorisation domain (UFD for
short) has been introduced and studied by Chatters and Jordan in \cite{ch,cj}.
Recently, the present authors, together with L Rigal, \cite{llr}, have shown
that many quantum algebras are noetherian UFD. In particular, we have shown
that the quantum group $O_q(SL_n)$ is a noetherian UFD.

Let $G$ be a connected simply connected complex semisimple
algebraic group. Since in the classical
setting it was shown by Popov, \cite{p}, that the ring of regular
functions on  $G$ is a unique factorisation domain, one can ask if a
similar result holds for the quantisation $O_q(G)$ of the coordinate
ring of $G$. The aim of this note is to provide a positive answer to
this question. In order to do this, we use a stratification of the prime
spectrum of $O_q(G)$ that was constructed by Joseph, \cite{josephbook}.

\section{Quantised enveloping algebras and quantum coordinate
  rings}\label{section:ufr}

\subsection{Quantised enveloping algebras}

Let $\g$ be a complex semisimple Lie algebra of rank $n$. We denote by
$\pi=\{\alpha_1,\dots,\alpha_n\}$ the set of simple roots associated to a
triangular decomposition $\g=\mathfrak{n}^- \oplus \mathfrak{h} \oplus
\mathfrak{n}^+$. Recall that $\pi$ is a basis of a euclidean vector space $E$
over $\mathbb{R}$, whose inner product is denoted by $(\mbox{ },\mbox{ })$
($E$ is usually denoted by $ \mathfrak{h}_{\mathbb{R}}^*$ in Bourbaki). We
denote by $W$ the Weyl group of $\g$; that is, the subgroup of the orthogonal
group of $E$ generated by the reflections $s_i:=s_{\alpha_i}$, for $i \in
\{1,\dots,n\}$, with reflecting hyperplanes $H_i:=\{\beta \in E \mid
(\beta,\alpha_i)=0\}$, for $i \in \{1,\dots,n\}$. If $w \in W$, we denote by
$l(w)$ its length. Further, we denote by $w_0$ the longest element of $W$.
Throughout this paper, the Coxeter group $W$ will be endowed with the Bruhat
order that we denote by $\leq$. We refer the reader to \cite[Appendix
A1]{josephbook} for the definition and properties of the Bruhat order.

We denote by  $R^+$ the set of positive roots and by $R$ the set of
 roots. We set $Q^+:=\mathbb{N} \alpha_1 \oplus \dots \oplus
 \mathbb{N} \alpha_n$. We denote by $\varpi_1,\dots,\varpi_n$ the fundamental
 weights, by $P$ the $\mathbb{Z}$-lattice generated by 
$\varpi_1,\dots,\varpi_n$, and by $P^+$ the set of dominant weights. 
In the sequel, $P$ will always be endowed with the following partial
 order:
$$ \lambda \leq \mu \mbox{ if and only if } \mu - \lambda \in Q^+.$$
Finally, we denote by $A=(a_{ij}) \in M_n(\mathbb{Z})$ the Cartan
 matrix associated to these data. 

Recall that the scalar product of two roots $(\alpha,\beta)$ is always
an integer. As in \cite{bg}, we assume that the short roots have
length $\sqrt{2}$.

For each $i \in \{1,\dots,n \}$, set
$q_i:=q^{\frac{(\alpha_i,\alpha_i)}{2}}$ and 
$$\left[ \begin{array}{l} m \\ k \end{array} \right]_i:=
\frac{(q_i-q_i^{-1}) \dots
  (q_i^{m-1}-q_i^{1-m})(q_i^m-q_i^{-m})}{(q_i-q_i^{-1})\dots
  (q_i^k-q_i^{-k})(q_i-q_i^{-1})\dots (q_i^{m-k}-q_i^{k-m})} $$
for all integers $0 \leq  k \leq  m$. By convention, we have  
$$\left[ \begin{array}{l} m \\ 0 \end{array} \right]_i:=1.$$

We will use the definition of the quantised enveloping algebra given in
\cite[I.6.3, I.6.4]{bg}. 
The quantised enveloping algebra $U_q(\g)$ of $\g$ over $\comp$ associated to
the previous data is the $\comp$-algebra generated by 
indeterminates $E_1,\dots,E_n,F_1,\dots , F_n,K_1^{\pm 1}, \dots, K_n^{\pm 1}$ subject to the following relations:
$$K_i K_j =K_j K_i \hspace{1cm} K_iK_i^{-1}=1$$
$$ K_i E_j K_i^{-1}=q_i^{a_{ij}}E_j  \hspace{1cm}  K_i F_j K_i^{-1}=q_i^{-a_{ij}}F_j$$
$$E_i F_j -F_jE_i=\delta_{ij} \frac{K_i-K_i^{-1}}{q_i-q_i^{-1}} $$
and the quantum Serre relations:
$$\sum_{k=0}^{1-a_{ij}} (-1)^k  \left[ \begin{array}{c} 1-a_{ij} \\ k
 \end{array} \right]_i E_i^{1-a_{ij} -k} E_j E_i^k=0  \mbox{ } (i \neq  j)$$
and 
$$\sum_{k=0}^{1-a_{ij}} (-1)^k  \left[ \begin{array}{c} 1-a_{ij} \\ k
 \end{array} \right]_i F_i^{1-a_{ij} -k} F_j F_i^k=0  \mbox{ } (i \neq  j).$$
Note that $U_q(\g)$ is a  Hopf algebra; its comultiplication is
defined by
$$\Delta(K_i)=K_i \otimes K_i \hspace{1cm} \Delta(E_i)=E_i \otimes 1 +
K_i\otimes E_i \hspace{1cm} \Delta(F_i)=F_i \otimes K_i^{-1} +
1 \otimes F_i, $$
its counit by 
 $$\varepsilon(K_i)= 1 \hspace{1cm}
 \varepsilon(E_i)=\varepsilon(F_i)=0, $$
and its antipode by 
$$S(K_i)=K_i^{-1} \hspace{1cm} S(E_i)=-K_i^{-1}E_i \hspace{1cm} S(F_i)=-F_iK_i.$$

We refer the reader to \cite{bg,jantzen,josephbook} for more details on
this algebra. 
Further, as usual, we denote by $U_q^+(\g)$ 
the
subalgebra of $U_q(\g)$ generated by $E_1,\dots,E_n$
and by  $U_q(\bk^+)$ 
the
subalgebra of $U_q(\g)$ generated by $E_1,\dots,E_n,K_1^{\pm 1},
\dots, K_n^{\pm 1}$.  
In a similar manner, $U_q^-(\g)$ is the subalgebra of $U_q(\g)$ generated by 
$F_1,\dots,F_n$ and $U_q(\bk^-)$ is the subalgebra of $U_q(\g)$ generated by
$F_1,\dots,F_n,K_1^{\pm 1},
\dots, K_n^{\pm 1}$. 

\subsection{Representation theory of quantised enveloping algebras}

It is well-known that the representation theory of the quantised
enveloping algebra $U_q(\g)$ is analogous to the representation theory of
the classical enveloping algebra $U(\g)$. In this section, we collect
the properties that will be needed in the rest of the paper.

As usual, if $M$ is a left $U_q(\g)$-module, we denote its dual by
$M^*$. Observe that $M^*$ is a right $U_q(\g)$-module in a natural
way. However, by using the antipode of $U_q(\g)$, this right
action of $U_q(\g)$ on $M^*$ can be twisted to a left action, so that
$M^*$ can be viewed as a left $U_q(\g)$-module.

Let $M$ be a $U_q(\g)$-module and $m \in M$. The element $m$ is said to have
weight $\lambda \in P$ if $K_i.m=q^{(\lambda,\alpha_i)} m$ for all $i \in
\{1,\dots,n\}$. For each $\lambda \in P$, set $$M_{\lambda}:=\{m\in M \mid
K_i.m=q^{(\lambda,\alpha_i)} m \mbox{ for all }i \in \{1,\dots,n\} \}.$$ If
$M_{\lambda} \neq 0$ then $M_{\lambda}$ is said to be a weight space of $M$
and $\lambda$ is a weight of $M$.

 It is well-known, see, for example \cite{bg,jantzen}, that, for each dominant
 weight $\lambda \in P^+$, there exists a unique (up to isomorphism) simple
 finite dimensional $U_q(\g)$-module of highest weight $\lambda$ that we
 denote by $V(\lambda)$. In the following proposition, we collect some
 well-known properties of the $V(\lambda)$, for $\lambda \in P^+$. We refer
 the reader to \cite[especially I.6.12]{bg}, \cite{humphreys} and
 \cite{jantzen} for details and proofs.

\begin{proposition}\label{irreductible}
Denote by
  $\Omega(\lambda)$ the set of those weights $\mu \in P$ such that
  $\vl_{\mu} \neq 0$.
\begin{enumerate}

\item $\vl = \oplus_{\mu \in \Omega(\lambda)}\,\vl_{\mu}$ 
\item The weights of $\vl$ are given by Weyl's character formula. 
In particular, if $\mu \in \Omega(\lambda)$, then $w\mu \in
\Omega(\lambda)$ for all $w \in W$.
\item For all $w \in W$, one has $\dime_{\comp} \vl_{w \lambda} = 1$.
\item $\vl^* \simeq V(-w_0 \lambda)$.
\item The weight $w_0 \lambda$ is the unique lowest weight of
  $\vl$. \\In particular, for all $\mu \in \Omega(\lambda)$, one has 
$w_0 \lambda \leq \mu \leq \lambda$.
\item $\Omega(\lambda)= \{\lambda - w \mu \mid w \in W \mbox{ and }
  \mu \in P^+ \mbox{ such that } \mu \leq \lambda \}$.

\end{enumerate}
\end{proposition}

For all $w \in W$ and $\lambda \in P^+$, let $u_{w\lambda}$ denote a
nonzero vector of weight $w\lambda$ in $\vl$. Then we denote by
$V_w^+(\lambda)$ the Demazure module associated to the pair $\lambda ,
w$, that is:
$$V_w^+(\lambda):= U_q^+(\g) u_{w\lambda}=U_q(\bk^+) u_{w\lambda}.$$
We also set 
$$V_w^-(\lambda):= U_q^-(\g) u_{w\lambda}=U_q(\bk^-) u_{w\lambda}.$$
(Observe that these definitions are independent of the choice of 
$u_{w\lambda}$ because of Proposition \ref{irreductible} (3).)

The following result may be well-known; however, we have been unable to locate
a precise statement.

\begin{proposition}\label{demazure}
\begin{enumerate}
\item $V_{w_0}^+(\lambda)=\vl=V_{id}^-(\lambda)$.
\item For all $i,j \in \{1,\dots,n\}$, one has 
$$V_{ w_0 s_i}^+(\varpi_j)= \left\{ \begin{array}{ll} 
\bigoplus_{\mu \in \Omega(\varpi_j) \setminus \{
  w_0\varpi_j \} } V(\varpi_j)_{\mu} & \mbox{ if } i=j \\
V(\varpi_j) & \mbox{ otherwise,}
\end{array}
\right.
$$
and 
$$V_{ s_i}^-(\varpi_j)= \left\{ \begin{array}{ll} 
\bigoplus_{\mu \in \Omega(\varpi_j) \setminus \{
  \varpi_j \} } V(\varpi_j)_{\mu} & \mbox{ if } i=j \\
V(\varpi_j) & \mbox{ otherwise.}
\end{array}
\right.
$$
\end{enumerate}
\end{proposition}
\begin{proof} We only prove the assertions corresponding to
``positive'' Demazure modules, the proof for ``negative'' Demazure modules is
similar.

Since $w_0 \lambda$ is the lowest weight of $\vl$, we have 
$U_q^+(\g) u_{w_0\lambda}=\vl$; that is, $V_{w_0}^+(\lambda)=\vl$. This
proves the first assertion.

In order to prove the second claim, we distinguish between two cases. 

First, let $i,j \in \{1,\dots,n\}$ with $i \neq j$. Then
$s_i(\varpi_j)=\varpi_j$. Hence, in this case, one has: $V_{ w_0
s_i}^+(\varpi_j)=U_q^+(\g) u_{ w_0 s_i\varpi_j}= U_q^+(\g) u_{ w_0
\varpi_j}=V_{ w_0 }^+(\varpi_j)= V(\varpi_j)$.

Next, let $j \in \{1,\dots,n\}$. 
Then $s_j (\varpi_j)=\varpi_j -
\alpha_j$. 
Let  $\mu \in \Omega(\varpi_j)$ with $\mu \neq  w_0\varpi_j$, and let 
$m \in V(\varpi_j)_{\mu}$ be any nonzero element. 
It follows from the first assertion
that there exists $x \in U_q^+(\g)$ such that $m=x. u_{w_0
  \varpi_j}$. The element $x$ 
  can be written as a linear combination of products 
$E_{i_1} \dots E_{i_k}$, with $k \in \mathbb{N}^*$ and $i_1,\dots,i_k
\in \{1,\dots,n\}$. Naturally, one can assume that $E_{i_1} \dots
E_{i_k}.u_{w_0 \varpi_j} \neq 0$ for each such product. 
Let $E_{i_1} \dots E_{i_k}$ be one of these 
products. 
Since $w_0 \pi =-\pi$, there exists $l \in \{1,\dots,n\}$ such that
$w_0 \alpha_{i_k}=-\alpha_l$. We will prove that $l=j$. Indeed,
assume that $l \neq j$. Since $E_{i_k}.u_{w_0 \varpi_j}$ is a
nonzero vector of $V(\varpi_j)$ of weight $w_0 \varpi_j +\alpha_{i_k}$,
we get that 
$$w_0 \varpi_j +\alpha_{i_k} \in \Omega(\varpi_j).$$
Then, we deduce from Proposition \ref{irreductible} that 
$$s_l w_0 \left(w_0 \varpi_j +\alpha_{i_k} \right) \in
\Omega(\varpi_j),$$
that is, 
$$s_l \varpi_j + \alpha_{l} \in \Omega(\varpi_j).$$
 Further, since we have assumed that $l \neq j$, we get
 $s_l\varpi_j=\varpi_j$, so that 
$$ \varpi_j + \alpha_{l} \in \Omega(\varpi_j).$$
This contradicts the fact that $\varpi_j$ is the highest weight of
$V(\varpi_j)$.

Thus, we have just proved that $w_0 \alpha_{i_k}=-\alpha_j$ for all
products $E_{i_1} \dots E_{i_k}$ that appear in $x$. 
Now, observe that $E_{i_k}.u_{w_0 \varpi_j}$ is a nonzero vector of
$V(\varpi_j)$ of weight $w_0 \varpi_j + \alpha_{i_k}=w_0 (\varpi_j
+w_0 \alpha_{i_k})= w_0 (\varpi_j -\alpha_j)=w_0 s_j \varpi_j$. 
Since $\dime_{\comp} V(\varpi_j)_{w_0s_j \varpi_j}=1$, we get that 
$E_{i_k}.u_{w_0 \varpi_j} = a u_{w_0s_j \varpi_j}$ for a certain
nonzero complex number $a$. Hence we get that 
$$m= x. u_{w_0  \varpi_j}=\sum \bullet E_{i_1} \dots E_{i_k}.u_{w_0 \varpi_j} =y.u_{w_0s_j \varpi_j},$$
where $\bullet$ denote some nonzero complex numbers and $y \in
U_q^+(\g)$. Thus $m \in V_{w_0s_j }^+(\varpi_j)$. This shows that 
$$\bigoplus_{\mu \in \Omega(\varpi_j) \setminus \{
  w_0\varpi_j \} } V(\varpi_j)_{\mu} \subseteq V_{ w_0
  s_j}^+(\varpi_j).$$
As the reverse inclusion is trivial, this finishes the proof.
\end{proof}

\subsection{Quantised coordinate rings of semisimple groups and their
  prime spectra.}

Let $G$ be a connected, simply connected, semisimple algebraic group
over $\comp$ with Lie algebra $\mathrm{Lie}(G)=\g$. Since $U_q(\g)$ is a Hopf algebra, one can define its Hopf dual 
$U_q(\g)^*$ (see \cite[1.4]{josephbook}) via
$$U_q(\g)^*:= \{ f \in \mathrm{Hom}_{\comp}(U_q(\g),\comp) \mid f=0
\mbox{ on some ideal of finite codimension} \}.$$

The quantised coordinate ring $O_q(G)$ of $G$ is the subalgebra of $U_q(\g)^*$
generated by the coordinate functions $c_{\xi, v}^{\lambda}$ for all $\lambda
\in P^+$, $\xi \in V(\lambda)^*$ and $v \in V(\lambda)$, where $c_{\xi,
v}^{\lambda}$ is the element of $U_q(\g)^*$ defined by $$c_{\xi,
v}^{\lambda}(u):= \xi (uv) \mbox{ for all } u \in U_q(\g),$$ see, for example,
\cite[Chapter 9]{josephbook}. As usual, if $\xi \in V(\lambda)_{\eta}^*$ and $v  \in
V(\lambda)_{\mu}$, we write $c_{\eta, \mu}^{\lambda}$ instead of $c_{\xi,
  v}^{\lambda}$. Naturally, this leads to some ambiguity. However,
when $\mu \in W.\lambda$ and $\eta \in W.(-w_0 \lambda)$, then 
$\dim (V(\lambda)_{\mu})=1=\dim(V(\lambda)_{\eta}^*)$, so that this
ambiguity is very minor.

It is well-known that $O_q(G)$ is a noetherian domain and a Hopf-subalgebra of $U_q(\g)^*$, see
\cite{bg,josephbook}. This latter structure allows us to define the so-called left and right
winding automorphisms (see, for instance, \cite[1.9.25]{bg} or
\cite[1.3.5]{josephbook}), and then to obtain an action of the torus
$\mathcal{H}:=(\comp^*)^{2n}$ on $O_q(G)$ (see \cite[5.2]{bgtrans}). 
More precisely, observe that the torus $H:=(\comp^*)^n$ can be
identified with $\mathrm{Hom}(P,\comp^*)$ via:
$$h(\lambda)=h_1^{\lambda_1} \dots h_n^{\lambda_n},$$
where $h=(h_1,\dots,h_n) \in H$ and $\lambda = \lambda_1 \varpi_1 +
\dots + \lambda_n \varpi_n$ with $\lambda_1, \dots, \lambda_n \in
\mathbb{Z}$. Then, it is known (see \cite[3.3]{hlt} or \cite[I.1.18]{bg}) 
that the torus
$\mathcal{H}$ acts rationally by $\comp$-algebra automorphisms on
$O_q(G)$ via: 
$$g.c_{\xi,  v}^{\lambda}= g_1(\mu) g_2(\eta) c_{\xi,
  v}^{\lambda},$$
for all $g=(g_1,g_2) \in \mathcal{H}=H \times H$, $\lambda \in
  P^+$, $\xi \in V(\lambda)_{\mu}^*$ and $v \in V(\lambda)_{\eta}$.
\\(We refer the reader to \cite[II.2.6]{bg} for the definition of a rational action.)

As usual, we denote by $\spec (O_q(G))$ the set of prime ideals in
$O_q(G)$. Recall that Joseph has proved \cite{josephjalgebra} that every prime in
$O_q(G)$ is completely prime.

Since $\mathcal{H}$ acts by automorphisms on $O_q(G)$, this induces an
action of $\mathcal{H}$ on the prime spectrum of $O_q(G)$. As usual,
we denote by $\mathcal{H}$-$\mathrm{Spec}(O_q(G))$ the set of those
primes ideals of $O_q(G)$ that are $\mathcal{H}$-invariant. This is a
finite set since Brown and Goodearl \cite[Section 5]{bgtrans} (see also
\cite[II.4]{bg}) have shown using previous results
of Joseph that 
$$ \ch \mbox{-}\spec (O_q(G))=\left\{ Q_{w_+,w_-} \mid (w_+,w_-) \in W
\times W \right\},$$
where 
$$Q_{w_+}^+:= \ideal{c_{\xi , v}^{\lambda} \mid \lambda \in P^+
  \mbox{, } v \in V(\lambda)_{\lambda} \mbox{ and } \xi \in
  (V_{w_+}^+(\lambda))^{\perp} \subseteq V(\lambda)^*   }, $$
$$Q_{w_-}^-:= \ideal{c_{\xi , v}^{\lambda} \mid \lambda \in P^+
  \mbox{, } v \in V(\lambda)_{w_0\lambda} \mbox{ and } \xi \in
  (V_{w_-w_0}^-(\lambda))^{\perp} \subseteq V(\lambda)^*   }, $$
and 
$$Q_{w_+,w_-}:=Q_{w_+}^+ + Q_{w_-}^-.$$
Since $q$ is transcendental, it follows from \cite[Th\'eor\`eme
  3]{josephcras} that it is enough to consider the fundamental weights
in the definition of $Q_{w_+}^+$ and $Q_{w_-}^-$. More precisely, we
deduce from  \cite[Th\'eor\`eme 3]{josephcras} the following result.

\begin{theorem}[Joseph]
\label{joseph}
$$ \ch \mbox{-}\spec (O_q(G))=\left\{ Q_{w_+,w_-} \mid (w_+,w_-) \in W
\times W \right\},$$
where 
$$Q_{w_+}^+:= \ideal{c_{\xi , v}^{\varpi_j} \mid j \in \{1, \dots, n\}
  \mbox{, } v \in V(\varpi_j)_{\varpi_j} \mbox{ and } \xi \in
  (V_{w_+}^+(\varpi_j))^{\perp} \subseteq V(\varpi_j)^*   }, $$
$$Q_{w_-}^-:= \ideal{c_{\xi , v}^{\varpi_j} \mid  j \in \{1, \dots, n\}
  \mbox{, } v \in V(\varpi_j)_{w_0 \varpi_j} \mbox{ and } \xi \in
  (V_{w_-w_0}^-(\varpi_j))^{\perp} \subseteq V(\varpi_j)^*   }, $$
and 
$$Q_{w_+,w_-}:=Q_{w_+}^+ + Q_{w_-}^-.$$
Moreover the prime ideals $Q_{w_+,w_-}$, for $(w_+,w_-) \in W \times W$,
are pairwise distinct.
\end{theorem}

\section{$O_q(G)$ is a noetherian UFD.}

In this section, we prove that $O_q(G)$ is a noetherian UFD (We refer the
reader to \cite[Section 1]{llr} for the definition of a noetherian UFD; the
key point is that each height one prime ideal should be generated by a normal
element.) In order to do this, we proceed in three steps.
\begin{enumerate}
\item First, by using results of Joseph, we show that there exist a finite
  number of nonzero normal $\ch$-eigenvectors $r_1,\dots,r_k$ of $O_q(G)$
  such that each $\ideal{r_i}$ is (completely) prime, and that 
  each nonzero
  $\ch$-invariant prime ideal of $O_q(G)$ contains one of the $r_i$. 
This property may be thought of as a ``weak factoriality'' result: 
$O_q(G)$ is
  an $\ch$-UFD in the terminology of \cite{llr}. 
\item Secondly, by 
using the $H$-stratification theory of Goodearl and Letzter (see
  \cite[II]{bg}), we show that the localisation of $O_q(G)$ with
  respect to the multiplicative system generated by the $r_i$ is a
  noetherian UFD.
\item Finally, we use a noncommutative analogue of Nagata's Lemma (see
  \cite[Proposition 1.6]{llr}) 
  to prove that $O_q(G)$ itself is a noetherian UFD. 
\end{enumerate}

\subsection{$O_q(G)$ is an $\ch$-UFD}

This aim of this section is two-fold. First, we show that 
for each $i \in \{ 1, \dots n\}$, the ideal generated by the normal
element 
$c_{-\varpi_i, w_0 \varpi_i}^{\varpi_i}$ 
or $c_{-w_0\varpi_i, \varpi_i}^{\varpi_i}$ is (completely) prime 
and then we prove that every
nonzero $\ch$-invariant prime ideal of $O_q(G)$ contains either one of
the $c_{-\varpi_i, w_0 \varpi_i}^{\varpi_i}$ or one of the
$c_{-w_0\varpi_i, \varpi_i}^{\varpi_i}$.

\begin{lemma}
\label{generatorheightoneHprimes}
Let $i \in \{ 1, \dots n\}$. 
Then $Q_{w_0,s_iw_0}= \ideal{c_{-\varpi_i ,
    w_0 \varpi_i}^{\varpi_i}}$ and $Q_{w_0s_i,w_0}= \ideal{c_{-w_0\varpi_i ,
     \varpi_i}^{\varpi_i}}$.
\end{lemma}
\begin{proof}
Recall that 
$$Q_{w_0,s_iw_0}= Q_{w_0}^++Q_{s_iw_0}^-,$$
where $$Q_{w_0}^+= \ideal{c_{\xi , v}^{\varpi_j} \mid j \in \{1, \dots, n\}
  \mbox{, } v \in V(\varpi_j)_{\varpi_j} \mbox{ and } \xi \in
  (V_{w_0}^+(\varpi_j))^{\perp} \subseteq V(\varpi_j)^*   }, $$
$$Q_{s_iw_0}^-= \ideal{c_{\xi , v}^{\varpi_j} \mid  j \in \{1, \dots, n\}
  \mbox{, } v \in V(\varpi_j)_{w_0 \varpi_j} \mbox{ and } \xi \in
  (V_{s_i}^-(\varpi_j))^{\perp} \subseteq V(\varpi_j)^*   }. $$
Next, it follows from Proposition \ref{demazure}(1) that 
$V_{w_0}^+(\varpi_j)=V(\varpi_j)$ for all $j$, so that 
$Q_{w_0}^+=(0)$. Also, we deduce from Proposition \ref{demazure}(2) that 
$V_{s_i}^-(\varpi_j)=V(\varpi_j)$ if $ j \neq i$, and 
   $V_{s_i}^-(\varpi_i)= \oplus_{\mu \in \Omega(\varpi_i)
  \setminus\{\varpi_i\} } V(\varpi_i)_{\mu}$. Hence, 
$$Q_{s_iw_0}^-= \ideal{c_{\xi , v}^{\varpi_i} 
\mid  v \in V(\varpi_i)_{w_0 \varpi_i} \mbox{ and } \xi \in
  V(\varpi_i)_{-\varpi_i}^*   }, $$
that is, $Q_{s_iw_0}^-= \ideal{c_{-\varpi_i , w_0 \varpi_i}^{\varpi_i}
  }$. Therefore 
  $ Q_{w_0,s_iw_0}= Q_{w_0}^++Q_{s_iw_0}^-=\ideal{c_{-\varpi_i ,
    w_0 \varpi_i}^{\varpi_i} }$, as desired.

The second claim of the lemma is obtained in the same way.
\end{proof}

Now observe that, in \cite{josephbook}, Joseph uses slighty different
 conventions for the dual $M^*$ of a left $U_q(\g)$-module. Indeed, it
 is mentioned in  \cite[9.1]{josephbook} that the dual $M^*$ is viewed
 with its natural right $U_q(\g)$-module structure. As a consequence,
 Joseph's convention for the weights of the dual $L(\lambda)^*$ of
 $L(\lambda)$, for $\lambda \in P^+$, is not exactly the same as our
 convention. In particular, the elements $c_{\varpi_i, w_0 \varpi_i}^{\varpi_i}$ and 
$c_{w_0\varpi_i, \varpi_i}^{\varpi_i}$, $i \in \{ 1, \dots ,n
\}$, that appear in  \cite[Corollary 9.1.4]{josephbook}, correspond 
to the elements $c_{-\varpi_i, w_0 \varpi_i}^{\varpi_i}$ and 
$c_{-w_0\varpi_i, \varpi_i}^{\varpi_i}$ in our notation. With this in mind, it follows from \cite[Corollary 9.1.4]{josephbook}
that the elements $c_{-\varpi_i, w_0 \varpi_i}^{\varpi_i}$ and 
$c_{-w_0\varpi_i, \varpi_i}^{\varpi_i}$, for $i \in \{ 1, \dots ,n
\}$, 
are normal in $O_q(G)$. Thus  we deduce
from Lemma \ref{generatorheightoneHprimes} the
following result which will allow us later to use a noncommutative
analogue of Nagata's Lemma in
order to prove that $O_q(G)$ is a noetherian UFD.

\begin{corollary}
\label{Hfacto}
The $2n$ elements $c_{-\varpi_i, w_0 \varpi_i}^{\varpi_i}$ and 
$c_{-w_0\varpi_i, \varpi_i}^{\varpi_i}$, for $i \in \{ 1, \dots ,n \}$, 
are
nonzero normal elements of $O_q(G)$ and they generate pairwise
distinct completely prime ideals of $O_q(G)$. 
\end{corollary}

Since the $c_{-\varpi_i, w_0 \varpi_i}^{\varpi_i}$ and $c_{-w_0\varpi_i,
\varpi_i}^{\varpi_i}$, for $i \in \{ 1, \dots ,n \}$, are $\ch$-eigenvectors
of $O_q(G)$, in order to prove that $O_q(G)$ is an $\ch$-UFD in the sense of
\cite[Definition 2.7]{llr}, it only remains to prove that
every nonzero $\ch$-invariant prime ideal of $O_q(G)$ contains either one of 
the  $c_{-\varpi_i, w_0 \varpi_i}^{\varpi_i}$ or one of the 
 $c_{-w_0\varpi_i, \varpi_i}^{\varpi_i}$.  This is
what we do next. 

\begin{lemma}
\label{inclusion}
Let $\mathbf{w}=(w_+,w_-) \in W \times W$, with $\mathbf{w} \neq
(w_0,w_0)$. Then $Q_{\mathbf{w}}$ contains either one of the 
$c_{-\varpi_i, w_0 \varpi_i}^{\varpi_i}$, or one of the
$c_{-w_0\varpi_i, \varpi_i}^{\varpi_i}$.
\end{lemma}

\begin{proof}
Since  $\mathbf{w} \neq (w_0,w_0)$, either $w_+ \neq w_0$, or $w_-
\neq w_0$. Assume, for instance, that $w_+ \neq w_0$, so that 
there exists $i \in \{ 1, \dots, n\}$ such that $w_+ \leq w_0s_i$. 
One can easily check from the definition of $Q_{\mathbf{w}}$ 
that this forces  $c_{-w_0\varpi_i,
  \varpi_i}^{\varpi_i} \in Q_{w_+}^+$, so that 
$$c_{-w_0\varpi_i,  \varpi_i}^{\varpi_i} \in Q_{w_+}^+ \subseteq Q_{\mathbf{w}},$$
as required.
\end{proof}

As a consequence of Corollary \ref{Hfacto} and Lemma \ref{inclusion}, we
get the following result.

\begin{corollary}
\label{HUFD}
$O_q(G)$ is an $\ch$-UFD.
\end{corollary}

\begin{proof}
Theorem~\ref{joseph} establishes that $ \ch \mbox{-}\spec (O_q(G))=\left\{
Q_{w_+,w_-} \mid (w_+,w_-) \in W \times W \right\}$. Note that $Q_{w_+,w_-} =
0$ precisely when $w_+ = w_- = w_0$. Thus, Corollary \ref{Hfacto} and Lemma
\ref{inclusion} show that each nonzero $\ch$-prime ideal of $O_q(G)$ contains
a nonzero $\ch$-prime of height one that is generated by a normal
$\ch$-eigenvector. Thus, $O_q(G)$ is an $\ch$-UFD.
\end{proof}

\subsection{$O_q(G)$ is a noetherian UFD.}

Set $T$ to be the localisation of $O_q(G)$ with respect to the
multiplicatively closed set generated by the normal $\ch$-eigenvectors
$c_{-\varpi_i, w_0 \varpi_i}^{\varpi_i}$ and $c_{-w_0\varpi_i,
\varpi_i}^{\varpi_i}$, for $i \in \{ 1, \dots ,n \}$. 
Then the rational action of
$\ch$ on $O_q(G)$ extends to an action of $\ch$ on the localisation $T$ by
$\comp$-algebra automorphisms, since we are localising with respect to
$\ch$-eigenvectors, and this action of $\ch$ on $T$ is also rational, by using
\cite[II.2.7]{bg}. The following result is a consequence of Corollary
\ref{HUFD} and \cite[Proposition 3.5]{llr}.

\begin{proposition}\label{t-is-hsimple}
The ring $T$ is $\ch$-simple; that
is, the only $\ch$-ideals of $T$ are $0$ and $T$. 
\end{proposition}

We are now in position to show that $O_q(G)$ is a noetherian UFD.

\begin{theorem}\label{oqgufd}
$O_q(G)$ is a noetherian UFD.  
\end{theorem}

\begin{proof} 
By \cite[Proposition 1.6]{llr}, it is enough to
prove that the localisation $T$ is a noetherian UFD. Now, as proved in
Proposition~\ref{t-is-hsimple}, $T$ is an $\ch$-simple ring. Thus, using
\cite[II.3.9]{bg}, $T$ is a noetherian UFD, as required.
\end{proof}

As a consequence, we deduce from Theorem \ref{oqgufd} and
\cite[Theorem 2.4]{cj} the following result.

\begin{corollary}
$O_q(G)$ is a maximal order.
\end{corollary}

The fact that $O_q(G)$ is a maximal order can also 
be proved directly by
using a suitable localisation of $O_q(G)$, 
\cite[Corollary 9.3.10]{josephbook}, which
is itself a maximal order.\\

\noindent{\bf Acknowledgment}~ We thank Laurent Rigal with whom we first
discussed this problem. We also thank Christian Ohn for a very helpful
conversation concerning the representation theory of $U(\g)$ during a meeting
of the Groupe de Travail Inter-universitaire en Alg\`ebre in La Rochelle and
thank the organisers for the opportunity to attend this meeting.



\vskip 1cm

\newpage

\noindent S Launois:\\
School of Mathematics, University of Edinburgh,\\
James Clerk Maxwell Building, King's Buildings, Mayfield Road,\\
Edinburgh EH9 3JZ, Scotland\\
E-mail : stephane.launois@ed.ac.uk
\\

\noindent T H Lenagan: \\
School of Mathematics, University of Edinburgh,\\
James Clerk Maxwell Building, King's Buildings, Mayfield Road,\\
Edinburgh EH9 3JZ, Scotland\\
E-mail: tom@maths.ed.ac.uk


\end{document}